\documentstyle[12pt,amscd,amsfonts]{amsart}

\newtheorem{thm}{Theorem}[section]
\newtheorem{cor}[thm]{Corollary}

\newtheorem{assa}{Assumption 1}

\newtheorem{assb}{Assumption 2}

\newtheorem{conj}{Conjecture}[section]
\newtheorem{mthm}{Main Theorem}

\newcommand{\DEG}{\operatorname{deg}}
\newcommand{\DEg}{\widehat{\DEG}}
\newcommand{\CH}{\widehat{CH}}
\newcommand{\FFF}{F_{\infty }}
\newcommand{\Mtrc}{\Vert \ \Vert}
\newcommand{\MTRC}{\Vert \ \Vert _{\sigma }}
\newcommand{\LLL}{L_{H,\Mtrc }}
\newcommand{\LLS}{L_{H,\MTRC }}
\newcommand{\IIm}{\operatorname{Im}}
\newcommand{\KER}{\operatorname{Ker}}
\newcommand{\fin}{\operatorname{fin}}
\newcommand{\COK}{\operatorname{Cok}}
\newcommand{\prim}{\operatorname{prim}}
\newcommand{\coprim}{\operatorname{coprim}}
\newcommand{\chern}{c_1(H, \Mtrc )}
\newcommand{\Chern}{\hat{c}_1(H, \Mtrc )}

\begin{document}

\title[Standard conjectures and their arithmetic analogues]
{A relation between standard conjectures and their arithmetic analogues}

\author{Yuichiro Takeda}

\address{Department of Mathematics, Faculty of science, 
Tokyo Metropolitan University, Minami-osawa 1-1, Hachioji-shi, 
Tokyo 192-03, Japan} 

\email{takeda@@math.metro-u.ac.jp}


\maketitle

\vskip 2pc
\section{Introduction}
\vskip 1pc

Let $X$ be a smooth projective variety of dimension $n$ defined 
over an algebraically closed field $k$.
We denote by $A^p(X)$ the vector space over $\Bbb R$ consisting 
of algebraic cycles of codimension $p$ modulo homological equivalence.
In \cite{grothendieck} Grothendieck conjectured that $A^p(X)$ 
behaves like complex cohomology.

\begin{conj}[Standard conjectures] \ \ \ \ 
Let $H$ be an ample line bundle on $X$ and 
$$
L_H:A^p(X)\to A^{p+1}(X)
$$
the homomorphism intersecting with the first Chern class $c_1(H)$.
Then for $p\leq \frac{n}{2}$, we have the followings:

\noindent
$\bold A_p(X, H)$: \ $L_H^{n-2p}:A^p(X)\to A^{n-p}(X)$ is an isomorphism.

\noindent
$\bold H_p(X, H)$: \ For $0\not= x\in A^p(X)$ such that $L_H^{n+1-2p}(x)=0$, 
$(-1)^p\DEG (L_H^{n-2p}(x)x)$ is positive.
\end{conj}

\vskip 1pc
$\bold A_p(X, H)$ is called the hard Lefschetz conjecture and 
$\bold H_p(X, H)$ is called the Hodge index conjecture.
When the characteristic of $k$ is zero, the Hodge index conjecture 
is already proved.

On the other hand, for an arithmetic variety the intersection theory 
of cycles was established by Arakelov \cite{arakelov} for surfaces and 
Gillet and Soul\'{e} \cite{GS1} for higher dimensional varieties.
It is quite natural to ask whether analogues of standard conjectures 
hold in this situation.
We now explain this.

Let $X$ be a regular scheme which is projective and flat over $\Bbb Z$.
We assume that the generic fiber $X_{\Bbb Q}$ is smooth over $\Bbb Q$.
Such a scheme is called an arithmetic variety.
For an arithmetic variety $X$ the arithmetic Chow group $\CH ^p(X)$ is 
defined and the intersection product on $\CH ^p(X)_{\Bbb Q}$ is 
established in \cite{GS1}.
We denote by $\FFF $ the complex conjugation on the complex manifold 
$X(\Bbb C)$ associated with the scheme $X\underset{\Bbb Z}{\otimes }\Bbb C$.
A line bundle $H$ on $X$ togather with an $\FFF $-invariant smooth 
hermitian metric $\Mtrc $ on the pull back $H_{\Bbb C}$ on $X(\Bbb C)$ 
is called a hermitian line bundle.
For a hermitian line bundle $(H, \Mtrc )$ on $X$, the arithmetic 
first Chern class $\Chern \in \CH ^1(X)$ is defined.
By intersecting with this class we obtain a homomorphism
$$
\LLL :\CH ^p(X)_{\Bbb R}\to \CH ^{p+1}(X)_{\Bbb R}.
$$
In \cite{GS} Gillet and Soul\'{e} proposed the following conjectures:

\begin{conj}[Arithmetic analogues of standard conjectures] \ \ \ \ 
Let $n$ be the relative dimension of $X$ over $\Bbb Z$ and $H$ an 
ample line bundle on $X$.
Then there exists an $\FFF $-invariant hermitian metric $\Mtrc $ on 
$H_{\Bbb C}$ satisfying the followings for $2p\leq n+1$:

\noindent
$\bold A_p(X, H, \Mtrc )$: \ $\LLL ^{n+1-2p}:\CH ^p(X)_{\Bbb R}\to 
\CH ^{n+1-p}(X)_{\Bbb R}$ is an isomorphism.

\noindent
$\bold H_p(X, H, \Mtrc )$: \ For $0\not= x\in \CH ^p(X)_{\Bbb R}$ such 
that $\LLL ^{n+2-2p}(x)=0$, $(-1)^p\DEg (\LLL ^{n+1-2p}(x)x)$ is positive.
\end{conj}

\vskip 1pc
When $n=1$, Conjecture 1.2 is derived from the Hodge index theorem for 
arithmetic surfaces by Faltings \cite{faltings} and Hriljac \cite{hriljac}.
This was extended to higher dimensional varieties by Moriwaki \cite{moriwaki}.
He proved $\bold H_1(X, H, \Mtrc )$ for any arithmetically ample 
hermitian line bundle $(H, \Mtrc )$ on an arithmetic variety $X$.
In \cite{kunnemann} K\"{u}nnemann showed that the conjectures are 
deduced from the similar conjectures for Arakelov Chow groups and 
proved them for projective spaces.

The aim of this paper is to resolve Conjecture 1.2 into other 
well-known conjectures including original standard conjectures.
In \cite{faltings, hriljac} for an arithmetic surface $X$ 
$\bold H_1(X, H, \Mtrc )$ is proved by the positivity of N\'{e}ron-Tate 
height pairing of the Jacobian of $X$ and by some arguments on the 
intersection of cycles of $X$ whose supports do not meet the generic 
fiber $X_{\Bbb Q}$.
In other words, our main results are generalisations of their methods 
to higher dimensional varieties.
As a consequence, all results as mentioned above are obtained, 
independent of the notion of the arithmetic ampleness.

\vskip 2pc
\section{Statements of the main results}
\vskip 1pc

We first recall some basic facts of Arakelov intersection theory.
Throughout the paper, $n$ is the relative dimension of an arithmetic 
variety $X$ over $\Bbb Z$.
For an arithmetic variety $X$, we put
\begin{gather*}
  Z^{p,p}(X)=\{\omega ; \text{real closed $(p,p)$-form on $X(\Bbb C)$ 
with $\FFF ^*\omega =(-1)^p\omega $}\}, \\
  H^{p,p}(X)=\{c\in H^{p,p}(X(\Bbb C)); \text{$c$ is real with 
$\FFF ^*c=(-1)^pc$}\}.
\end{gather*}
Then we have an exact sequence
\begin{align*}
CH^{p-1,p}(X)_{\Bbb R}\overset{\rho }{\to }&H^{p-1,p-1}(X)\overset{a}
{\to }\CH ^p(X)_{\Bbb R}\overset{(\zeta ,\omega )}{\longrightarrow } \\
&CH^p(X)_{\Bbb R}\oplus Z^{p,p}(X)\overset{cl-h}{\to }H^{p,p}(X)\to 0,
\end{align*}
where the definitions of $CH^{p,p-1}, a, \zeta $ and $\omega $ are 
seen in \cite[3.3]{GS1}.
The map $cl:CH^p(X)_{\Bbb R}\to H^{p,p}(X)$ is the cycle class map and 
$h$ is the canonical projection map.
The map $\rho :CH^{p-1,p}(X)_{\Bbb R}\to H^{p-1,p-1}(X)$ is the regulator 
map up to constant factor by \cite[Theorem 3.5]{GS1}.

We fix a smooth $\FFF $-invariant K\"{a}hler metric $h$ on $X(\Bbb C)$.
The pair $\overline{X}=(X, h)$ is called an Arakelov variety.
By identifying an element of $H^{p,p}(X)$ with a harmonic $(p,p)$-form 
with respect to $h$, we can regard $H^{p,p}(X)$ as a subspace of $Z^{p,p}(X)$.
We put $CH^p(\overline{X})=\omega ^{-1}(H^{p,p}(X))$ and call it Arakelov 
Chow group.

\begin{assa} \ \ \ \ 
For $p\geq 0$, the vector space $H^{p,p}(X)$ is spanned by images of $\rho $ 
and $cl$, that is,
$$
H^{p,p}(X)=\IIm \rho \oplus \IIm cl.
$$
\end{assa}

\vskip 1pc

This is a part of Beilinson conjectures.
If Assumption 1 holds for $X$ and for $p-1$, by the definition of the 
Arakelov Chow group we obtain the following exact sequence:
$$
CH^{p-1}(X)_{\Bbb R}\overset{a\cdot cl}{\to }CH^p(\overline{X})_{\Bbb R}
\overset{\zeta }{\to }CH^p(X)_{\Bbb R}\to 0.
$$

Let $(H, \Mtrc )$ be a hermitian line bundle on $X$.
Suppose that $H$ is ample and that $\Mtrc $ is a positive metric.
Then the first Chern form of $(H, \Mtrc )$ determines an $\FFF $-invariant 
K\"{a}hler metric $h$ on $X(\Bbb C)$.
Since the product with the K\"{a}hler form respects harmonicity of forms, 
for the Arakelov variety $\overline{X}=(X, h)$ we can define a homomorphism
$$
\LLL :CH^p(\overline{X})_{\Bbb R}\to CH^{p+1}(\overline{X})_{\Bbb R}.
$$

\begin{conj} \ \ \ \ 
For an ample line bundle $H$ on $X$, there exists a positive $\FFF $-
invariant hermitian metric $\Mtrc $ on $H_{\Bbb C}$ satisfying the 
followings for $2p\leq n+1$:

\noindent
$\bold A\bold A_p(X, H, \Mtrc )$: \ $\LLL ^{n+1-2p}: CH^p(\overline{X})
_{\Bbb R}\to CH^{n+1-p}(\overline{X})_{\Bbb R}$ is an isomorphism.

\noindent
$\bold A\bold H_p(X, H, \Mtrc )$: \ For $0\not= x\in CH^p(\overline{X})
_{\Bbb R}$ such that $\LLL ^{n+2-2p}(x)=0$, $(-1)^p\DEg (\LLL ^{n+1-2p}(x)x)$ 
is positive.
\end{conj}

\vskip 1pc
\begin{thm} \ \ \ \ 
$\bold A\bold A_p(X, H, \Mtrc )$ implies $\bold A_p(X, H, \Mtrc )$.
$\bold A\bold H_p(X, H, \Mtrc )$ implies $\bold H_p(X, H, \Mtrc )$.
\end{thm}

\vskip 1pc
Theorem 2.1 was proved by K\"{u}nnemann in \cite{kunnemann}.
He also proved Conjecture 2.1 for projective spaces.
Recently the author proved Conjecture 2.1 for regular quadric hypersurfaces 
in \cite{takeda}.

From now on, every arithmetic variety $X$ is assumed to be irreducible.
Then $X$ is defined over the ring of integers $\cal O_K$ of an algebraic 
number field $K$ and the generic fiber $X_K$ is geometrically irreducible.
The cycle class map $cl:CH^p(X)_{\Bbb R}\to H^{p,p}(X)$ factors through 
the Chow group of the generic fiber and the restriction map $CH^p(X)\to 
CH^p(X_K)$ is surjective.
Hence the image of $cl$ coincides with the image of $CH^p(X_K)_{\Bbb R}$.
We denote it by $A^p(X_K)$.
Then the preceding exact sequence yields
$$
0\to A^{p-1}(X_K)\to CH^p(\overline{X})_{\Bbb R}\overset{\zeta }{\to }
CH^p(X)_{\Bbb R}\to 0.
$$

For an ample line bundle $H$ on $X$, we can consider the hard Lefschetz 
conjecture and the Hodge index conjecture for $A^p(X_K)$ although $X_K$ 
is not defined over an algebraically closed field.
We denote them by $\bold A_p(X_K, H_K)$ and $\bold H_p(X_K, H_K)$ 
respectively.
For an algebraic closure $\overline{K}$ of $K$, standard conjectures 
for $(X_{\overline{K}}, H_{\overline{K}})$ imply these for $(X_K, H_K)$.
In particular, $\bold H_p(X_K, H_K)$ is true.

We put
$$
CH_{\fin }^p(X)=\KER (CH^p(X)\to CH^p(X_K)).
$$
For an ample line bundle $H$ on $X$, we can define a homomorphism
$$
L_H:CH_{\fin }^p(X)_{\Bbb R}\to CH_{\fin }^{p+1}(X)_{\Bbb R}
$$
by intersecting with the first Chern class $c_1(H)$.
For $x\in CH_{\fin }^p(X)_{\Bbb R}$ we denote by $\tilde{x}\in CH^p
(\overline{X})_{\Bbb R}$ a lifting of $x$.
Then we define a pairing
$$
\langle \ , \ \rangle :CH_{\fin }^p(X)_{\Bbb R}\otimes CH_{\fin }
^{n+1-p}(X)_{\Bbb R}\to \Bbb R
$$
by $\langle x, y\rangle =\DEg (\tilde{x}\tilde{y})$.
This definition is independent of the choice of liftings.
Here we propose the following conjectures:

\begin{conj} \ \ \ \ 

\noindent
$\bold F\bold A_p(X, H)$: \ $L_H^{n+1-2p}:CH_{\fin }^p(X)_{\Bbb R}\to 
CH_{\fin }^{n+1-p}(X)_{\Bbb R}$ is an isomorphism.

\noindent
$\bold F\bold H_p(X, H)$: \ For $0\not=x\in CH_{\fin }^p(X)_{\Bbb R}$ 
such that $L_H^{n+2-2p}(x)=0$, $(-1)^p\langle L_H^{n+1-2p}(x), x\rangle $ 
is positive.
\end{conj}

We now construct the height pairing by means of Arakelov intersection theory.
We put
$$
CH^p(X_K)^0_{\Bbb R}=\KER (cl:CH^p(X_K)_{\Bbb R}\to H^{p,p}(X)).
$$
We need the following hypothesis.

\begin{assb} \ \ \ \ 
For any $x\in CH^p(X_K)^0_{\Bbb R}$, there exists a lifting $\tilde{x}
\in CH^p(\overline{X})_{\Bbb R}$ such that $\tilde{x}$ is orthogonal 
to every vertical cycles.
That is to say, it holds that
$$
\DEg (\tilde{x}\cdot (C, 0))=0
$$
for any cycle $C$ of codimension $n+1-p$ whose support does not meet 
the generic fiber $X_K$.
\end{assb}

\vskip 1pc
If $CH^p(X_K)^0_{\Bbb R}$ and $CH^{n+1-p}(X_K)^0_{\Bbb R}$ admit the 
above assumption, then we can define a pairing
$$
\langle \ , \ \rangle :CH^p(X_K)^0_{\Bbb R}\otimes CH^{n+1-p}(X_K)^0
_{\Bbb R}\to \Bbb R
$$
by $\langle x, y\rangle =\DEg (\tilde{x} \tilde{y})$, where $\tilde{x}$ 
and $\tilde{y}$ are liftings which satisfy the above assumption.
This definition is independent of the choice of the liftings.

Suppose that the pairing
$$
\langle \ , \ \rangle :CH_{\fin }^p(X)_{\Bbb R}\otimes CH_{\fin }
^{n+1-p}(X)_{\Bbb R}\to \Bbb R
$$
is nondegenerate.
For a lifting $\tilde{x}\in CH^p(\overline{X})_{\Bbb R}$ of $x\in 
CH^p(X_K)^0_{\Bbb R}$, the homomorphism
$$
CH_{\fin }^{n+1-p}(X)_{\Bbb R}\to \Bbb R, \ y\mapsto \DEg (\tilde{x}\tilde{y})
$$
is determined independently of the choice of the lifting $\tilde{y}$.
Because of the nondegeneracy of the pairing there exists a unique 
$z\in CH_{\fin }^p(X)_{\Bbb R}$ such that
$$
\DEg (\tilde{x}\tilde{y})=\langle z, y\rangle 
$$
for any $y\in CH_{\fin }^{n+1-p}(X)_{\Bbb R}$.
Then for a lifting $\tilde{z}$, $\tilde{x}-\tilde{z}$ is a lifting of 
$x\in CH^p(X_K)^0_{\Bbb R}$ which satisfies Assumption 2.
In the same way we can show that $CH^{n+1-p}(X_K)^0_{\Bbb R}$ also admits 
Assumption 2 and we can define the height pairing.
In particular, if $\bold F\bold A_k(X, H)$ and $\bold F\bold H_k(X, H)$ 
hold for $k\leq p$, then Assumption 2 holds for $p$ and $n+1-p$.

\vskip 1pc
{\it Remark 1}: \ 
We assume that the pairing
$$
\langle L_H^{n+1-2p} \ \ , \ \rangle :CH_{\fin }^p(X)_{\Bbb R}\otimes 
CH_{\fin }^p(X)_{\Bbb R}\to \Bbb R
$$
is nondegenerate.
Then we can say that $CH^{n+1-p}(X_K)^0_{\Bbb R}$ admits Assumption 2 and 
that there exists a lifting $\tilde{x}\in CH^p(\overline{X})_{\Bbb R}$ of 
arbitrary $x\in CH^p(X_K)^0_{\Bbb R}$ such that $\LLL ^{n+1-2p}(\tilde{x})$ 
is orthogonal to any vertical cycles.
Therefore we can define the pairing
$$
\langle L_{H_K}^{n+1-2p} \ \ , \ \rangle :CH^p(X_K)^0_{\Bbb R}\otimes 
CH^p(X_K)^0_{\Bbb R}\to \Bbb R
$$
although $CH^p(X_K)^0_{\Bbb R}$ may not admit Assumption 2.

\vskip 1pc
For the height pairing for $X$, Beilinson conjectured an analogues of 
standard conjectures in \cite[\S 5]{beilinson}.

\begin{conj} \ \ \ \ 
For $2p\leq n+1$, we assume either Assumption 2 for $p$ and $n+1-p$ or 
nondegeneracy of the pairing of $CH_{\fin }^p(X)_{\Bbb R}$ in Remark 1.
Then we have the followings:

\noindent
$\bold H\bold A_p(X, H)$: \ $L_{H_K}^{n+1-2p}:CH^p(X_K)^0_{\Bbb R}\to 
CH^{n+1-p}(X_K)^0_{\Bbb R}$ is an isomorphism.

\noindent
$\bold H\bold H_p(X, H)$: \ For $0\not=x\in CH^p(X_K)^0_{\Bbb R}$ such that 
$L_{H_K}^{n+2-2p}(x)=0$, $(-1)^p\langle L_{H_K}^{n+1-2p}(x), x\rangle $ is 
positive.
\end{conj}

\vskip 1pc
{\it Remark 2}: \ 
We now choose another positive $\FFF $-invariant hermitian metric 
$\Mtrc ^{\prime }$ on $H_{\Bbb C}$.
Then there exists an $\FFF $-invariant positive real valued function $f$ 
on $X(\Bbb C)$ such that $\Mtrc ^{\prime }=f\Mtrc $.
We denote the K\"{a}hler metric associated with $\Mtrc ^{\prime }$ by 
$h^{\prime }$ and let $\overline{X}^{\prime }=(X, h^{\prime })$.
For $x\in CH_{\fin }^p(X)_{\Bbb R}$, we choose a lifting $\tilde{x}\in 
CH^p(\overline{X})_{\Bbb R}$ of $x$.
Since $\omega (\tilde{x})=0$, $\tilde{x}$ is also contained in the Arakelov 
Chow group $CH^p(\overline{X}^{\prime })_{\Bbb R}$ of $\overline{X}^{\prime }$.
Hence the definition of pairing of $CH_{\fin }^p(X)_{\Bbb R}$ does not 
depend on the choice of positive metrics on $H_{\Bbb C}$.
In the same way it can be shown that the definition of the pairing of 
$CH^p(X)^0_{\Bbb R}$ and Assumption 2 do not also depend on the choice of 
metrics.
Moreover since
\begin{align*}
 \DEg (\hat{c}_1(H, \Mtrc ^{\prime }&)^{n+1-2p}\tilde{x}^2)=\DEg 
((\Chern -2a(\log f))^{n+1-2p}\tilde{x}^2)  \\
 &=\DEg (\Chern ^{n+1-2p}\tilde{x}^2-2(n+1-2p)a(\log f\chern ^{n-2p}
\omega (\tilde{x})^2))  \\
 &=\DEg (\Chern ^{n+1-2p}\tilde{x}^2),
\end{align*}
the conjectures $\bold F\bold H_p(X, H)$ and $\bold H\bold H_p(X, H)$ 
do not depend on the choice of metrics.

\vskip 1pc
Here we state our main theorem.

\begin{mthm} \ \ \ \ 
Let $X$ be an arithmetic variety defined over the ring of integers 
$\cal O_K$ of an algebraic number field $K$.
Suppose that the generic fiber $X_K$ is geometrically irreducible.
Let $H$ be an ample line bundle on $X$.
Given a positive $\FFF $-invariant hermitian metric $\Mtrc $ on 
$H_{\Bbb C}$, we define a metric $\MTRC $ by $\MTRC =\exp (\sigma )\Mtrc $ 
for $\sigma \in \Bbb R$.

\noindent
i) \ We assume Assumption 1 for $p-1$ and $n-p$.
Then $\bold F\bold A_p(X, H)$, $\bold H\bold A_p(X, H)$, 
$\bold A_p(X_K, H_K)$ and $\bold A_{p-1}(X_K, H_K)$ imply 
$\bold A_p(X, H, \MTRC )$ for almost all $\sigma $.

\noindent
ii) \ We assume Assumption 1 for $p-1$.
We suppose that the pairings
$$
\langle L_H^{n+1-2p} \ \ , \ \rangle :CH_{\fin }^p(X)_{\Bbb R}
\otimes CH_{\fin }^p(X)_{\Bbb R}\to \Bbb R
$$
and
$$
\langle \ , \ \rangle :CH_{\fin }^{p-1}(X)_{\Bbb R}\otimes 
CH_{\fin }^{n+2-p}(X)_{\Bbb R}\to \Bbb R
$$
are nondegenerate.
Then $\bold F\bold H_p(X, H)$, $\bold H\bold H_p(X, H)$ and 
$\bold A_{p-1}(X_K, H_K)$ imply \linebreak
$\bold H_p(X, H, \MTRC )$ for $0 \ll -\sigma $.
\end{mthm}

\vskip 1pc
After arithmetic ampleness of a hermitian line bundle was defined in 
\cite{zhang}, a finer version of Conjecture 1.2 is proposed.
This says that for arithmetically ample hermitian line bundle $(H, \Mtrc )$ 
on $X$, $\bold A_p(X, H, \Mtrc )$ and $\bold H_p(X, H, \Mtrc )$ hold.
For arbitrary hermitian line bundle $(H, \Mtrc )$ which satisfies the 
conditions in Main Theorem and for $0\ll -\sigma $, the hermitian line 
bundle $(H, \MTRC )$ becomes arithmetically ample.
But it is difficult to compare the upper bound of $\sigma $ such that 
$(H, \MTRC )$ is arithmetically ample with that for which ii) of Main 
Theorem holds.

\begin{cor} \ \ \ \ 
If $X$ is a Grassmannian or a projective smooth toric scheme or a 
generalized flag scheme, which is a quotient scheme of a split reductive 
group scheme by a Borel subgroup, defined over a ring of integers $\cal O_K$.
Then Conjecture 1.2 holds for any ample line bundle $H$ with a positive 
metric $\Mtrc $ and for $0\ll -\sigma $.
\end{cor}

{\it Proof}: \ \ 
Since $X$ can be stratified into finite pieces isomorphic to affine spaces 
over $\cal O_K$, Assumption 1 holds.
Since $CH_{\fin }^p(X)_{\Bbb R}$ and $CH^p(X_K)^0_{\Bbb R}$ vanish for all 
$p$, all assumptions about these are vacuous.
For an embedding $\tau :K\to \Bbb C$, we denote by $X_{\tau }$ the complex 
manifold associated with the scheme $X\underset{\tau }{\otimes }\Bbb C$.
Then the cycle class map $cl:CH^p(X_K)_{\Bbb R}\to H^{p,p}(X_{\tau })$ is 
bijective.
So the standard conjectures for $X_K$ are valid.
Hence Main Theorem implies Conjecture 2.2.
\qed

\vskip 1pc
\begin{cor} \ \ \ \ 
Let $X$ be an arithmetic variety.
Then the Hodge index conjecture $\bold H_1(X, H, \MTRC )$ of codimension 
one holds for any ample line bundle $H$ with a positive metric $\Mtrc $ and 
for $0\ll -\sigma $.
\end{cor}

{\it Proof}: \ \ 
We have only to verify the conditions in Main Theorem for $p=1$.
We denote by $\Sigma _K$ the set of all infinite places of $K$.
Then $H^{0,0}(X)\simeq \underset{\sigma \in \Sigma _K}{\oplus }\Bbb R$ 
and $CH^{0,1}(X)\simeq \cal O_K^{\times }$.
Hence by Dirichlet unit theorem Assumption 1 for $p=0$ holds.
Nondegeneracy of the pairing $\langle L_H^{n-1} \ \ , \ \rangle $ on 
$CH_{\fin }^1(X)_{\Bbb R}$ and negativity of $\langle L_H^{n-1}x, x\rangle $ 
for $0\not= x\in CH_{\fin }^1(X)_{\Bbb R}$ have been already proved in 
\cite[Lemma 1.3]{moriwaki}.
Since $CH^0(X)_{\Bbb R}\simeq \Bbb R$ and $CH^{n+1}(X)_{\Bbb R}=0$, we have 
$CH_{\fin }^0(X)_{\Bbb R}=CH_{\fin }^{n+1}(X)_{\Bbb R}=0$.
$\bold H\bold H_1(X, H)$ holds by the positivity of N\'{e}ron-Tate height 
pairing and $\bold A_0(X_K, H_K)$ is trivial.
Hence the proof has completed.
\qed

\vskip 2pc
\section{Proof of the main theorem}
\vskip 1pc

By Theorem 2.1 we have only to prove $\bold A\bold A_p(X, H, \MTRC )$ and 
$\bold A\bold H_p(X, H, \MTRC )$.
We begin by the proof of i).
Since $\bold F\bold A_p(X, H)$, $\bold H\bold A_p(X, H)$ and 
$\bold A_{p-1}(X_K, H_K)$ hold, 
$$
L_H^{n+1-2p}:CH^p(X)_{\Bbb R}\to CH^{n+1-p}(X)_{\Bbb R}
$$
is surjective and its kernel is isomorphic to
$$
A_{\prim }^p(X_K)=\KER (L_{H_K}^{n+1-2p}:A^p(X_K)\to A^{n+1-p}(X_K)).
$$
We put
$$
A_{\coprim }^{n-p}(X_K)=\COK (L_{H_K}^{n+1-2p}:A^{p-1}(X_K)\to A^{n-p}(X_K))
$$
and consider the following diagram:
\begin{equation*}
 \begin{CD}
  @. 0 @. @. A_{\prim }^p(X_K) @.  \\
  @. @VVV @. @VVV @. \\
  0 @>>> A^{p-1}(X_K) @>>> CH^p(\overline{X})_{\Bbb R} @>>> 
CH^p(X)_{\Bbb R} @>>> 0  \\
  @. @VV{L_{H_K}^{n+1-2p}}V @VV{\LLS ^{n+1-2p}}V @VV{L_H^{n+1-2p}}V @. \\
  0 @>>> A^{n-p}(X_K) @>>> CH^{n+1-p}(\overline{X})_{\Bbb R} @>>> 
CH^{n+1-p}(X)_{\Bbb R} @>>> 0 \\
  @. @VVV @. @VVV @. \\
  @. A_{\coprim }^{n-p}(X_K) @. @. 0. @.
 \end{CD}
\end{equation*}
To show that $\LLS ^{n+1-2p}$ is bijective, we have only to prove that 
the edge homomorphism of the above diagram is an isomorphism.
Given $x\in A_{\prim }^p(X_K)$, we choose a lifting $\tilde{x}\in 
CH^p(\overline{X})_{\Bbb R}$ of $x$.
Then we have
\begin{align*}
 \LLS ^{n+1-2p}(\tilde{x})&=\hat{c}_1(H, \MTRC )^{n+1-2p}\tilde{x} \\
 &=(\Chern -2a(\sigma ))^{n+1-2p}\tilde{x}  \\
 &=(\Chern ^{n+1-2p}-2\sigma (n+1-2p)a(\chern ^{n-2p}))\tilde{x} \\
 &=\Chern ^{n+1-2p}\tilde{x}-2\sigma (n+1-2p)a(\chern ^{n-2p}x).
\end{align*}
Hence if we denote by $F_{\sigma }$ the edge homomorphism of the above 
diagram for $\sigma $, then we have
$$
F_{\sigma }=F_0-2\sigma (n+1-2p)L_H^{n-2p}.
$$
$\bold A_p(X_K, H_K)$ and $\bold A_{p-1}(X_K, H_K)$ imply that 
$L_H^{n-2p}:A_{\prim }^p(X_K)\to A_{\coprim }^{n-p}(X_K)$ is an isomorphism.
Hence for all but finitely many $\sigma $, $F_{\sigma }$ is an isomorphism.

We turn to the proof of ii).
Let $x$ be a primitive element in $CH^p(\overline{X})_{\Bbb R}$ with 
respect to $\LLL $, that is, $\LLL ^{n+2-2p}(x)=0$ holds.
Then for a closed $(p,p)$-form $\omega $ we have
\begin{align*}
 \LLS ^{n+2-2p}(x+a(\omega ))&=(\Chern -2a(\sigma ))^{n+2-2p}(x+a(\omega )) \\
 &=(\Chern ^{n+2-2p}-2\sigma (n+2-2p)a(\chern ^{n+1-2p}))(x+a(\omega )) \\
 &=a(\chern ^{n+2-2p}\omega -2\sigma (n+2-2p)\chern ^{n+1-2p}\omega (x)).
\end{align*}
Hence $x+a(\omega )\in CH^p(\overline{X})_{\Bbb R}$ is primitive with 
respect to $\LLS $ if and only if
$$
\chern ^{n+2-2p}\omega =2\sigma (n+2-2p)\chern ^{n+1-2p}\omega (x).
$$
Since $\chern ^{n+2-2p}\omega (x)=0$, $\bold A_{p-1}(X_K, H_K)$ implies 
that $\omega (x)$ is written by 
$$
\begin{cases}
\omega _0(x)+\chern \omega _1(x) &\text{if} \ 2p<n+1   \\
\chern \omega _1(x) &\text{if} \ 2p=n+1,
\end{cases}
$$
where each $\omega _i(x)$ is contained in $A_{\prim }^{p-i}(X_K)$.
If $2p=n+1$, then we set $\omega _0(x)=0$.
By the above equality we have
$$
\omega =2\sigma (n+2-2p)\omega _1(x).
$$
Hence we can say that the vector spaces consisting of primitive cycles in 
$CH^p(\overline{X})_{\Bbb R}$ with respect to $\LLS $ are isomorphic for 
every $\sigma \in \Bbb R$.

We first assume $\omega (x)=0$.
Then the restriction of $\zeta (x)\in CH^p(X)_{\Bbb R}$ to the generic fiber 
$X_K$ is homologically trivial.
We denote the restriction of $\zeta (x)$ by $x^{\prime }$.
Let $x_1$ be a lifting of $x^{\prime }$ which is orthogonal to $L_H^{n+1-2p}
CH_{\fin }^p(X)_{\Bbb R}$ and $x_0=x-x_1$.
We suppose that $x_0$ is represented by $(Z, 0)$ and the support of the 
cycle $Z$ does not meet the generic fiber $X_K$.
We identify $x_0$ with the cycle $[Z]\in CH_{\fin }^p(X)_{\Bbb R}$.

Since $\LLS ^{n+2-2p}(x)=0$, we have
$$
\LLS ^{n+2-2p}(x_1)=-L_H^{n+2-2p}(x_0)\in CH_{\fin }^{n+2-p}(X)_{\Bbb R}.
$$
Since $\langle \ , \ \rangle :CH_{\fin }^{p-1}(X)_{\Bbb R}\otimes 
CH_{\fin }^{n+2-p}(X)_{\Bbb R}\to \Bbb R$ is nondegenerate and 
$\LLS ^{n+2-2p}(x_1)$ is orthogonal to $CH_{\fin }^{p-1}(X)_{\Bbb R}$, 
the both sides of the above equality are zero.
Hence we have
\begin{align*}
 (-1)^p\DEg (\LLS ^{n+1-2p}(x)x)&=(-1)^p\langle L_H^{n+1-2p}(x_0), x_0
\rangle  \\
 & \ +(-1)^p\langle L_{H_K}^{n+1-2p}(x^{\prime }), x^{\prime }\rangle 
\end{align*}
and it is positive by $\bold F\bold H_p(X, H)$ and $\bold H\bold H_p(X, H)$.

We next consider when $\omega (x)\not= 0$.
Then any primitive cycle $y$ in $CH^p(\overline{X})_{\Bbb R}$ with respect 
to $\LLS $ is written by $x+a(\omega )$ where $x$ is a primitive cycle with 
respect to $\LLL $ and $\omega =2\sigma (n+2-2p)\omega _1(x)$.
Then we have
\begin{align*}
\LLS ^{n+1-2p}(y)y&=(\Chern -2a(\sigma ))^{n+1-2p}(x+a(\omega ))^2 \\
 &=(\Chern ^{n+1-2p}-2\sigma (n+1-2p)a(\chern ^{n-2p}))(x^2+2a(\omega (x)
\omega ))). 
\end{align*}
By substituting the above equality for $\omega $ and $\omega (x)=
\omega _0(x)+\chern \omega _1(x)$ we have
\begin{multline*}
\DEg (\LLS ^{n+1-2p}(y)y)=\DEg (\Chern ^{n+1-2p}x^2)-\sigma (n+1-2p)
\DEG (\chern ^{n-2p}\omega _0(x)^2)   \\
+\sigma (n+3-2p)\DEG (\chern ^{n+2-2p}\omega _1(x)^2).
\end{multline*}
The Hodge index theorem for $(X_K, H_K)$ implies
$$
(-1)^p\DEG (\chern ^{n-2p}\omega _0(x)^2)>0
$$
and
$$
(-1)^{p+1}\DEG (\chern ^{n+2-2p}\omega _1(x)^2)>0
$$
if $\omega _i(x)\not= 0$.
When $2p<n+1$, $\omega (x)\not= 0$ implies $\omega _0(x)\not= 0$ or 
$\omega _1(x)\not= 0$.
Hence for $0 \ll -\sigma $, $(-1)^p\DEg (\LLS ^{n+1-2p}(x)x)$ is positive.
When $2p=n+1$, $\omega _0(x)=0$ and $\omega _1(x)\not= 0$.
Hence in this case the same inequality is obtained.
\qed


\vskip 2pc





\begin{thebibliography}{9}





\bibitem{arakelov}S. J. Arakelov,
Intersection theory of divisors on an arithmetic surface, 
{\it Math. USSR Izv.} {\bf 8} (1974), 1167--1180.


\bibitem{beilinson}A. A. Beilinson,
Height pairing between algebraic cycles, in {\it K-Theory, Arithmetic 
and Geometry} (ed. Yu. I. Manin), Lect. Note Math. vol.1289, Springer, 
Berlin Heidelberg, 1987, pp. 1--26.


\bibitem{bloch}S. Bloch,
Height pairings for algebraic cycles, {\it J. Pure Appl. Algebra} 
{\bf 34} (1984), 119--145.


\bibitem{faltings}G. Faltings,
Calculus on arithmetic surfaces, {\it Ann. Math.} {\bf 119} (1991), 387--424.


\bibitem{GS1}H. Gillet and C. Soul\'{e},
Arithmetic intersection theory, {\it Publ. Math. I.H.E.S.} {\bf 72} 
(1990), 93--174.


\bibitem{GS}\bysame,
Arithmetic analogs of standard conjectures, in {\it Motives} (eds. 
U. Jannsen, S. Keiman and J.-P. Serre), Proc. Sym. Pure Math. vol.55 
part 1, American Mathematical Society, Providence, 1994, pp. 129--140.


\bibitem{grothendieck}A. Grothendieck,
Standard conjectures on algebraic cycles, in {\it Algebraic Geometry 
(Papers Presented at the Bombey Colloquium, 1968)}, Oxford University 
Press, Oxford, 1969, pp. 193--199.


\bibitem{hriljac}P. Hriljac,
Heights and Arakelov's intersection theory, {\it Amer. J. Math.} {\bf 107} 
(1985), 23--38.


\bibitem{kunnemann}K. K\"{u}nnemann,
Some remarks on the arithmetic Hodge index theorem, {\it Compositio Math.} 
{\bf 99} (1995), 109--128.


\bibitem{moriwaki}A. Moriwaki,
Hodge index theorem for arithmetic cycles of codimension one, {\it Math. 
Res. Lett.} {\bf 3} (1996), 173--183.


\bibitem{takeda}Y. Takeda,
Arithmetic analogues of standard conjectures for quadric hypersurfaces, 
preprint.


\bibitem{zhang}S. Zhang,
Positive line bundles on arithmetic varieties, {\it J. Amer. Math. Soc.} 
{\bf 8} (1995), 187--221.



\end{thebibliography}
\end{document}